\documentclass[12pt]{article}
\usepackage{amssymb,a4}
\usepackage{amsthm}
\usepackage{amsfonts}
\usepackage{mathrsfs}
\usepackage{amsmath}

\def\z{\mathbb{Z}}
\def\Z{\mathbb{Z}}
\def\c{\mathbb{C}}
\def\C{\mathbb{C}}
\def\dim{\hbox{dim}}

\def\ll{\lambda}

\newfont{\df}{eufm10}

\def\ll{\lambda}

\def\ot{\otimes}

\def\de{\delta}
\def\dim{\hbox{\rm dim}\,}

\def\Vir{\hbox{\rm Vir}}

\def\Vir{\hbox{Vir}}

\def\ot{\otimes}

\hoffset \voffset \oddsidemargin=60pt \evensidemargin=45pt
\numberwithin{equation}{section}

\title{\bf  Classification of Harish-Chandra
modules over some Lie algebras related to the Virasoro algebra
\footnote{E-mail:liudong@hutc.zj.cn }}
\author{Dong Liu\\ Department of Mathematics, Huzhou University\\ Zhejiang Huzhou, 313000, China
}
\date{ }

\begin{document}
\maketitle

\begin{abstract}
In this paper, we provide a uniform method to thoroughly classify
all Harish-Chandra modules over some Lie algebras related to the
Virasoro algebras. We first classify such modules over the Lie
algebra $W(\varrho)[s]$ for $s=0,\frac12$. With this result and
method, we can also do such works for some Lie algebras related to the Virasoro algebra, including
the several kinds of Schr\"{o}dinger-Virasoro Lie algebras, which
are open up to now.
\end{abstract}



\smallskip\bigskip

\newtheorem{theo}{Theorem}[section]
\newtheorem{defi}[theo]{Definition}
\newtheorem{lemm}[theo]{Lemma}
\newtheorem{coro}[theo]{Corollary}
\newtheorem{prop}[theo]{Proposition}

\section{Introduction}
\setcounter{section}{1}\setcounter{theo}{0} \setcounter{equation}{0}

Recently many infinite dimensional Lie algebras related to the
Virasoro algebra were studied sufficiently. Among them, the twisted
Heisenberg-Virasoro algebra was first studied by Arbarello et al. in
\cite{ADKP}, where a connection is established between the second
cohomology of certain moduli spaces of curves and the second
cohomology of the Lie algebra of differential operators of order at
most one. The $W$-algebra $W(2, 2)$ was introduced in \cite{ZD} for
the study of classification of vertex operator algebras generated by
weight 2 vectors. Schr\"{o}dinger-Virasoro algebras, playing
important roles in mathematics and statistical physics, were first
introduced by M. Henkle in \cite{H} by studying the free
Schr\"{o}dinger equation.

Since all the above Lie algebras are closely related to the Virasoro
Lie algebra, it is highly expected that their Hraish-Chandra modules
structures are well classified as the Virasoro algebra in \cite{M} (also see \cite{Kac}, \cite{MP}, \cite{MZ}, \cite{S}, etc.), the high rank Virasoro algebra in \cite{LvZ0} and
\cite{S1}, the Weyl algebra in \cite{S2}. However, up to now, there
are few results about it. In \cite{LvZ} all
Harish-Chandra modules over the twisted Heisenberg-Virasoro algebra
were classified. However their calculations are
very complicated and cannot to be used in general. In \cite{LS},
such study was considered for the original and twisted
Schr\"{o}dinger-Virasoro Lie algebras. However all irreducible
uniform bounded modules over them are yet not classified, it is the
key point to classify all Harish-Chandra modules and is open up to
now for many Lie algebras. So it is a very important question that how to
classify all Harish-Chandra modules, especially all irreducible
uniform bounded modules, over some Lie algebras related to the
Virasoro algebra. In this paper, we provide a new and uniform method
to thoroughly classify all Harish-Chandra modules over some Lie
algebras related to the Virasoro algebras, including all the above
Lie algebras. Throughout the paper, we shall use $\c, {\mathbb Q}$, $\z$ and $\z_+$ to denote the sets of complex numbers, rational
numbers, integers and positive integers, respectively. For any set $S$, we use
$S^*$ to denote the set of nonzero elements in $S$.

First we introduce the following Lie algebras.

\begin{defi} For $s=0, \frac12$ and $\varrho\in\mathbb Q$, as a vector space over $\c$, the Lie algebra $\mathcal L[s]=W(\varrho)[s]$ has a
basis
$\{L_n, Y_p\mid n \in\z,\; p\in\z+s\}$
with the following relations
\begin{eqnarray}
&&[L_m,L_n]=(n-m)L_{m+n},\\
&&[L_m,Y_p]=(p-m\varrho)Y_{m+p}, \\
&&[Y_p,Y_q]=0
\end{eqnarray} for all $m,\; n\in\z$ and $p,\; q\in\z+s$.
\end{defi}

The Lie algebra $W(\varrho)[s]$ can be realized from the semi-product of
the centerless Virasoro algebra $\Vir$ and the $\Vir$-module ${\cal F}_{\varrho}$ of
the intermediate series. In fact, let
Vir$=\hbox{Span}_\c\{L_m\mid m\in\z\}$ be the centerless Virasoro algebra (also called Witt algebra) and $H=\c\{Y_p\mid p\in\z+s\}$ (denote by ${\cal F}_{\varrho}$ in \cite{OR})  be a
$\Vir$-module with actions $L_m\cdot Y_p=(p-m\varrho)Y_{m+p}$ for any $m\in\z, \;p\in\z+s$, then $W(\varrho)[s]$ is just the Lie algebra
$\Vir\ltimes H$.

The Lie algebra ${\mathcal L}[s]$ is more connected to the the cohomology of the Virasoro algebra and extensions of some Lie algebras (see \cite{Fuk}, \cite{OR},\cite{OR2} and \cite{GJP}).

As Vir-modules, ${\cal F}_{0}$ and ${\cal F}_{-1}$ have $\c v_0$ as a submodule and quotient module respectively. Their corresponding irreducible quotient and submodule are isomorphic (both isomorphic to ${\mathcal A}_{0, 0}'$, see Section 2 in detail), so we always suppose that $\varrho\ne -1$ throughout this paper. Moreover if $s=\frac12$, we concentrate our study on the case of $\rho=\frac{m}{n}$, $(m, n)=1$ and $n$ is an even number, which is enough for almost Lie algebras as we known.

 The Lie algebra
$W(0)[0]$ is just the centerless twisted Heisenberg-Virasoro algebra (see \cite{ADKP}). It can be realized as the Lie algebra of differential operators of order at most 1.
The structure and representation theory of this Lie algebra and its universal central extensions were studied in many papers (see \cite{ADKP, B, LJ, LvZ, SJ, LWZ, SS, LPZ}, etc.).

The universal central extension of $W(1)[0]$ (named $W(2,2)$ in \cite{ZD}) can be also realized from the so-called
{\it loop-Virasoro algebra} $\mathcal G$ (see \cite{GLZ}). Let $\c[t, t^{-1}]$
be the Laurents polynomial ring over $\c$, then the $W(2,2)={\mathcal G}/(t^2)$.

The universal central extension of $W(1)[0]$ and its highest weight modules enter the picture
naturally during the discussion on $L(\frac12,0)\otimes L(\frac12,0)$. Its highest weight modules produce a new class of
vertex  operator algebras. Contrast to the Virasoro algebra case,
this class of vertex operator algebras are always irrational. Several papers studied its representation theory (see \cite{ZD, LGZ, JP, GLZ}, etc.).

Schr\"{o}dinger-Virasoro algebras (see Section 5 for their definitions) are more connected to $W(\varrho)[s]$. In fact, $W(0)[0]$ (resp. $W(\frac12)[s]$) is a subalgebra (resp. quotient algebra) of the original or twisted Schr\"{o}dinger-Virasoro algebra. Recently many researches on the  structure and representation theory of the original, twisted, deformative Schr\"{o}dinger-Virasoro algebras (see \cite{H}, \cite{DDM}, \cite{HU1}, \cite{HU2}, \cite{RU}, \cite{LS}, \cite{TZ}, \cite{WL}, \cite{DS}, etc.).

Therefore Lie algebras ${\mathcal L}[s]$ play very important roles
on the study of Lie algebras related to the Virasoro algebra. The
present paper is devoted to determining all Harish-Chandra modules
(i.e. irreducible weight modules with finite dimensional weight
spaces) over Lie algebras $W(\varrho)[s]$. More precisely we prove
that there are three different classes of Harish-Chandra modules
over them. One class is formed by simple modules of intermediate
series, the other two classes consist of the highest, lowest weight
modules. It is consistent with that for the Virasoro algebra case in
\cite{M}. Here we get a key lemma (Lemma 3.1 in Section 3) for the irreducible uniform bounded weight modules and then obtain the main result over the Lie algebra $W(\varrho)[s]$. With this method and result, we get a beautiful
application to determine irreducible uniform bounded weight modules
over Schr\"{o}dinger-Virasoro algebras.

The paper is arranged as follows. In Section 2, we recall some
notations and collect known facts about irreducible, indecomposable
modules over the Virasoro algebra. In Section 3, we determine all
irreducible uniformly bounded weight modules over ${{\mathcal
L}[s]}$ which turn out to be modules of intermediate series.
Moreover we obtain the main result of this paper: the classification
of irreducible weight ${{\mathcal L}[s]}$-modules with finite
dimensional weight space. As we mentioned, they are irreducible
highest, lowest weight modules, or irreducible modules of the
intermediate series. In Section 4, using the methods and results in
Section 3, we classified all Harish-Chandra modules over some Lie
algebra related the Virasoro algebra, including the original,
twisted and deformative Schr\"{o}dinger-Virasoro algebras.

For convenient, all modules considered in this paper are nontrivial. We always denote by $U(L)$ the universal enveloping algebra of a given Lie algebra $L$.

\section{Basics}

In this section, we collect some known facts for later use.

Introduce a $(1-s)\z$-gradation on ${\mathcal L}[s]={W(\varrho)[s]}$ by defining the degrees: deg $L_n=n$, deg $Y_p=p$.
Set
$${\mathcal L}[s]_+=\sum_{n, p>0}(\c L_n+\c Y_p),\,\,{\mathcal L}[s]_-=\sum_{n, p<0}(\c L_n+\c Y_p),$$
and
$${\mathcal L}[s]_{0}=\c L_0+\c(1-2s)Y_0.$$

The $(1-s)\z$-gradation of ${\mathcal L}[s]=\oplus_{p\in(1-s)\z} {\mathcal L}[s]_p$ induces $(1-s)\z$-gradations on the universal enveloping algebra $U({\mathcal L}[s])=\oplus_{p\in(1-s)\z} U({\mathcal L}[s])_p$ and the universal enveloping algebra $U(H)=\oplus_{p\in(1-s)\z} U(H)_p$, where $H$ is the abelian subalgebra of ${\mathcal L}[s]$ spanned by
$Y_r$ for all $r\in\z+s$.

The universal central extensions of ${\mathcal L}[0]$ were given in \cite{GJP}. Such work for ${\mathcal L}[\frac12]$ can also be easily done as \cite{GJP}.
\begin{prop}\cite{GJP}
$$H^2(W(\varrho)[0], \c)=\begin{cases}\c\gamma_0+\c\gamma_{01}+\c\gamma_{02},\ \varrho=0;\\ \c\gamma_0+\c\gamma_{11},\ \varrho=1;\\ \c\gamma_0, \ \hbox{otherwise}\end{cases}$$
with $2$-cocycles as follows:

$\gamma_{0}(L_m, L_n)=\frac1{12}(m^3-m)\delta_{m+n, 0}$;

$\gamma_{01}(L_m, I_n)=(m^2-m)\de_{m+n, 0}$;  $\gamma_{02}(I_m, I_n)=n\de_{m+n, 0}$;

$\gamma_{11}(L_m, I_n)=\frac1{12}(m^3-m)\delta_{m+n, 0}$.

\end{prop}

\vskip 5pt An ${\mathcal L}[s]$-module $V$ is called a weight module if
$V$ is the sum of all its weight spaces $V^\lambda=\{v\in V\mid L_0v=\lambda v\}$. For a weight module $V$
we define
$$\hbox{Supp}(V):=\bigl\{\lambda\in \c \bigm|V^\lambda\neq
0\bigr\},$$ which is called the weight set (or the
support) of $V$.

\par
 An irreducible weight ${\mathcal L}[s]$-module $V$ is called the intermediate
 series if all its weight space are one dimensional.
\par
A  weight ${\mathcal L}[s]$-module $V$ is called a {highest} (resp.
{lowest) weight module} with {highest weight} (resp. {lowest
weight}) $\lambda\in \c$, if there exists a nonzero weight vector $v
\in V^\lambda$ such that

1) $V$ is generated by $v$ as  ${\mathcal L}[s]$-module;

2) ${\mathcal L}[s]_+ v=0 $ (resp. ${\mathcal L}[s]_- v=0 $).

Obviously, if $V$ is an irreducible weight ${\mathcal L}[s]$-module, then
there exists $\lambda\in\c$ such that ${\rm
Supp}(V)\subset\lambda+(1-s)\z$. So $V=\sum_{p\in (1-s)\z}V_p$ is a $(1-s)\z$-graded module, where $V_p=V^{\ll+p}$.

If, in addition, all weight spaces $V^\ll$ of a weight ${\mathcal L}[s]$-module $V$ are finite dimensional, the module $V$ is called a {\it
Harish-Chandra module}. Clearly a highest (lowest) weight module
is a Harish-Chandra module.

Kaplansky-Santharoubane \cite{KS} in 1983 gave a classification of
$Vir$-modules of the intermediate series. There are three families
of indecomposable modules with each weight space is one-dimensional:

(1) ${\mathcal A}_{a,\; b}=\sum_{i\in\z}\c v_i$:
$L_mv_i=(a+i+b m)v_{m+i}$;

(2) ${\mathcal A}(a)=\sum_{i\in\z}\c v_i$: $L_mv_i=(i+m)v_{m+i}$
if $i\ne 0$, $L_mv_0=m(m+a)v_{m}$;

(3) ${\mathcal B}(a)=\sum_{i\in\z}\c v_i$:  $L_mv_i=iv_{m+i}$ if
$i\ne -m$, $L_mv_{-m}=-m(m+a)v_0$,  for some $a, b\in\c$.

It is well-known that ${\cal A}_{a,\; b}\cong{\cal A}_{a+1,\; b}, \forall a, b\in\c$, then we can always suppose that $a\not\in\z$ or $a=0$ in ${\cal A}_{a,\; b}$.
Moreover the module
${\cal A}_{a,\; b}$ is simple if $a\notin\z$ or $b\ne0, 1$.
 In the opposite case the module
contains two simple subquotients namely the trivial module and
$\c[t, t^{-1}]/\c$. It is also clear that ${\mathcal A}_{0,0}$ and
${\mathcal B}(a)$ both have $\c v_0$ as a submodule, and their corresponding
quotients are isomorphic, which we denote by ${\mathcal A}_{0,0}'$. Dually,
${\mathcal A}_{0,1}$ and ${\mathcal A}(a)$ both have $\C v_0$ as a quotient module, and
their corresponding submodules are isomorphic to ${\mathcal A}_{0,0}'$. For
convenience we simply write ${\mathcal A}_{a,b}'={\mathcal A}_{a,b}$ when ${\mathcal A}_{a,b}$ is
irreducible.

All Harish-Chandra modules over the Virasoro algebra were classified in \cite{M} (also in \cite{MP} and \cite{S}). Since then such works were done on the twisted Heisenberg-Virasoro algebra in \cite{LvZ}.

\begin{theo}(\cite{M}, etc)
Let $V$ be an irreducible weight Vir-module with finite dimensional weight spaces. Then $V$ is a highest weight module, lowest weight module, or Harish-Chandra module of intermediate series.
\end{theo}

\begin{theo} \cite{LvZ}
 Let $V$ be an irreducible weight module over
$W(0)[0]$  with all weight spaces finite-dimensional. Then $V$ is a highest weight module, a
lowest weight module, or a Harish-Chandra module of intermediate series.
\end{theo}

\noindent{\bf Remark.}
For $a, b, c\in\c$, let ${\mathcal A}_{a,\; b,\; c}$ be the $W(0)[0]$-module induced from the
Vir-module ${\mathcal A}_{a,\; b}$ with $Y_pv_k=cv_{p+k}$ for all $p, k\in\z$. Clearly  ${\mathcal A}_{a,\; b,\; c}$ is irreducible if and only if $a\notin\z$ or $b\ne 0, 1$ or $c\ne 0$. We also use ${\mathcal A}_{a,\; b,\; c}'$ to denote by the nontrivial simple subquotient of ${\mathcal A}_{a,\; b,\; c}$ $($or itself if it is irreducible$)$. From \cite{LvZ} or \cite{LJ} we know that any Harish-Chandra $W(0)[0]$-module of intermediate series is isomorphic to ${\mathcal A}_{a,\; b,\; c}'$ for some $a, b, c\in\c$.

In \cite{LS},
Harish-Chandra modules over the original and twisted
Schr\"{o}dinger-Virasoro Lie algebras were studied. However all irreducible
uniform bounded modules over them are yet not classified, it is the
key point to classify all Harish-Chandra modules and is open up to
now.

\begin{theo} \cite{LS}
 Let $V$ be an irreducible weight module over the Schr\"{o}dinger-Virasoro Lie algebra $\frak{sv}[s]$
 with all weight spaces finite-dimensional. Then $V$ is a highest weight module, a
lowest weight module, or a uniformly bounded irreducible weight module.
\end{theo}

\section{Harish-Chandra modules over ${\mathcal L}[s]$}

In this section,  we shall classify all irreducible weight module with finite dimensional weight spaces over ${\mathcal L}[s]$.

\begin{lemm} \label{lemm31}
Let $V$ be a uniformly bounded
nontrivial irreducible weight module over ${\mathcal L}[s]$. Suppose that there exist a nonzero $u_{i_0}\in V_{i_0}$ and $Y_{r_0},\; r_0\ne 0$ such that $Y_{r_0}u_{i_0}=0$. Then $HV=0$, where $H$ be the abelian subalgebra of ${\mathcal L}[s]$ spanned by
$Y_r$ for all $r\in\z+s$.\end{lemm}
\noindent{\bf Proof.}

Assume that $V=\sum V_p$ is a uniformly bounded
nontrivial irreducible weight module over ${\mathcal L}[s]$, where Supp$(V)=\{a+(1-s)\z\}$ and $V_p=\{v\in V\mid L_0v_p=(a+p)v\}$ for some $a\in\c$.

From
representation theory of $\Vir$ (see \cite{KS}), $\dim
V_p=n$ for all $a+p \neq 0$. If $a \in \z$, we can assume
 that $a=0$. Moreover, as a $\Vir$-module, $V$ has a $\Vir$-submodule
filtration
    $$0=W^{(0)}\subset W^{(1)} \subset W^{(2)}\subset \cdots \subset W^{(p)}=V,$$
where $W^{(1)}, \cdots ,W^{(p)}$ are $\Vir$-submodules of $V$, and
the quotient modules \break $W^{(i)}/W^{(i-1)}\cong {\mathcal A}_{a_i, b_i}'$ for some $a_i, b_i\in\C$.

Set $U=U(\Vir)u_{i_0}$. Clearly $U$ is a nontrivial Vir-module (if not, $U$ is a trivial ${\mathcal L}[s]$-module).
So $U$ is also a uniform bounded Vir-module and then there exists an irreducible Vir-submodule $V'\cong {\mathcal A}_{a, b}'$ of $U$ for some $a, b\in\c$ by the above statement and Theorem 2.2.  Moreover $V'=\sum_{i\in\z} \c v_i$ and $V=U(H)V'=\sum_{i\in\z} U(H)v_i$ (If $V'\cong {\mathcal A}_{0,0}'$, then $V=U(H)V'=\sum_{i\in\z^*} U(H)v_i$).

Fixed $i\in\z$ and $v_i\in U(\Vir)u_{i_0}$, there exists $n\in\z_+$ such that for every $v_i=f_iu_{i_0}$, where $f_i$ is a finite sum of terms $c_iL_{m_1}L_{m_2}\cdots L_{m_k}$ with all $k\le n-1$ and $c_i\in\c$.

Since $Y_{r_0}u_{i_0}=0$ then $Y_{r_0}U(H)u_{i_0}=0$. Moreover by $Y_{r_0}L_mu_{i_0}=L_mY_{r_0}u_{i_0}-(r_0-m\varrho)Y_{m+r_0}u_{i_0}$ we have $Y_{r_0}^2L_mu_{i_0}=0$.
Similarly $Y_{r_0}^{k+1}L_{m_1}L_{m_2}\cdots L_{m_k}u_{i_0}=0$.

So we have $Y_{r_0}^nv_i=0$.  Now since $L_jv_i=(a+bj+i)v_{i+j}$ for any $j\in\z$, then there exists $N\in\z_+$ such that $Y_{r_0}^NV'=0$.

 Therefore
\begin{equation}\label{eq31}Y_{r_0}^NV=0.\end{equation}
By actions of suitable $L_m$ for some $m\in\z$ on (\ref{eq31}) we can deduce that \begin{equation}\label{eq32}H^NV=0,\end{equation} where $H^N:=\c\{Y_{r_1}Y_{r_2}\cdots Y_{r_N}\mid r_i\in\z+s, i=1,2,\cdots, N\}$ (just as a vector space over $\c$).  Since $HV$ is a ${\mathcal L}[s]$-submodule of $V$, then $HV=V$ or $HV=0$. In the first case we also have $HV=0$ by (\ref{eq32}). \qed

\begin{theo}\label{T41}  Let $V$ be a uniformly bounded irreducible
${\cal L}[s]=W(\varrho)[s]$-module.
If $\varrho\ne 0$, then $V\simeq {\cal A}_{a,\;b,\;0}'$ for some
$a, b\in\c$, where ${\mathcal A}_{a,\; b,\; 0}$ is the ${\cal L}[0]$-module induced from the
Vir-module ${\mathcal A}_{a,\; b}$ with $Y_pv_k=0$ for all $p, k\in\z$, and
${\mathcal A}_{a,\; b,\; 0}'$ is nontrivial simple subquotient of ${\mathcal A}_{a,\; b,\; 0}$ $($or itself if it is irreducible$)$. Consequently, $V$ is an irreducible weight module over
the Virasoro algebra.

\end{theo}

\noindent{\bf Proof.}  We will use the notations from Lemma \ref{lemm31}.

By Lemma \ref{lemm31} we only need to consider the case that all $Y_i,\ i\in\z+s, i\ne 0$, are nonsingular on $V$.
In the case of $s=0$, if there exists a nonzero $v\in V$ such that $Y_0v=0$, then as in the proof of Lemma \ref{lemm31}, we can also deduce that $H^nV=0$. It is contradict to the hypothetical conditions. So $Y_0$ is also nonsingular on $V$.

So we can suppose that all $Y_i, i\in\z+s$, are nonsingular on $V$. Then all $V_p=\{v\in V\mid L_0v=(a+p)v\}, p\in(1-s)\z$, have a same finite dimension.

Now choose an irreducible Vir-submodule $V'=\sum\c  v_i\cong \mathcal A_{a, b}'$ of $V$ and then $V=U(H)V'=\sum U(H)v_i$.

Let $h=Y_{i}Y_{-i}\in U(H)_0$ and $k\in(1-s)\z$, then $h$ is a linear transformation over the finite dimensional vector space $V_k$.
So there exists a minimal polynomial $m(x)$ of $h$ such that $m(h)V_k=0$ (if $s={1\over 2}$, then we can regard $h$ as a linear transformation over the finite dimensional vector space $V_k+V_{k+\frac12}$). By actions of some $Y_j$ we can get
\begin{equation}m(h)V=0.\label{nipl}\end{equation}

Since $\varrho\in\mathbb Q^*$, so there exist $m_0\in\z$ and $i_0\in\z+s$ such that $i_0-m_0\varrho=0$ (where we suppose that if $s=\frac12$, then $\rho=\frac{m}{n}$, $(m, n)=1$ and $n$ is an even number, see Section 1).
By action of $L_{m_0}$  on (\ref{nipl}) for $h=Y_{i_0}Y_{-i_0}$, we have $Y_{i_0}Y_{m_0-i_0}m'(h)V=0$ then $m'(h)V=0$. It gets a contradiction. \qed

\noindent{\bf Remark.} The above method is not suitable for the case of $\varrho=0$.

\begin{prop}\label{P43} Let  V be an irreducible weight module over
$W(\varrho)[s]\, (\varrho\ne0)$ with all weight spaces finite-dimensional. If V is not
uniformly bounded, then V is either a highest weight module or a
lowest weight module.
\end{prop}
\noindent{Proof.} The proof is essentially same as that in the case of $\varrho=0$ in \cite{LvZ} (also see \cite{LS} for the $\varrho=\frac12$ case).
\qed

 With Theorems \ref{T41} and Proposition \ref{P43}, we obtain the main result of this paper.

\begin{theo}\label{Main} If $V$ is an irreducible weight
module over $W(\varrho)[s] \,(\varrho\ne 0)$, with finite dimensional weight spaces, then $V$
is a
highest or lowest weight module or the Harish-Chandra module of the intermediate series.
\end{theo}

\section{Harish-Chandra modules over some Lie algebras}

In this section, we use the results and methods in Section 3 to classify all Harish-Chandra modules over some Lie algebras related to the Virasoro algebra. \emph{Especially, we get a beautiful application to the Schr\"odinger algebras for such work, which was conjectured several years.}

\subsection{The Schr\"{o}dinger-Virasoro Lie algebra}

For $s=0, \frac12$, the Schr\"{o}dinger-Virasoro Lie algebra $\frak{sv}[s]$
introduced in \cite{HU1,HU2,RU}, in the context of
non-equilibrium statistical physics as a by-product of the
computation of $n$-point functions that are covariant under the
action of the Schr\"{o}dinger group, is the infinite-dimensional Lie
algebra with $\c$-basis $\{L_n,\,M_n,\,Y_p\mid n\in
\z,\, p\in\z+s\}$ and Lie brackets,
\begin{eqnarray}
\!\!\!&\!\!\!& [L_m,\,L_{n}]=(n-m)L_{m+n} ,\label{1eq1.1}
\\
\!\!\!&\!\!\!&[L_m,\,Y_p]=(p-\frac{m}{2})Y_{p+m} ,\label{1eq1.2}\\[2pt]
\!\!\!&\!\!\!& [L_m,\,M_n]=nM_{n+m} ,\label{1eq1.3}\\[2pt]
\!\!\!&\!\!\!&[Y_p,\,Y_{q}]=(q-p)M_{q+p} ,\label{1eq1.4}\\[2pt]
\!\!\!&\!\!\!&[Y_p,\,M_n]=[M_n,\,M_{n'}]=0.\label{1eq1.5}
\end{eqnarray}

The Lie algebras $\frak{sv}[\frac12]$ and $\frak{sv}[0]$ are called the original and twisted Schr\"{o}dinger-Virasoro algebra respectively in \cite{RU}.
We also denote by $H=\c\{Y_p, p\in\z+s\}$.

Let $V$ be an irreducible weight ${\frak{sv}[s]}$-module with finite dimensional weight spaces. If $V$ is not a highest (lowest) weight module, then it is uniformly bounded (see \cite{LS}). However, it is the most difficulty in classifying Harish-Chandra modules to classify all irreducible uniformly bounded modules. It is opened up to now. With the methods and results in Section 3, we can easily do such works.

\begin{prop} \label{p51}
Let $V$ be a uniformly bounded irreducible
$\frak{sv}[s]$-module. Then $M_nV=Y_{n+s}V=0$ for all $n\in\z$.

\end{prop}

\noindent{\bf Proof.}
Set $\mathcal T=\c\{L_m, M_n\mid n\in\z\}$ is the subalgebra of $\frak{sv}[s]$, then $\mathcal T\cong W(0)[0]$.
Choose an irreducible $\mathcal T$-submodule $V'\cong{\mathcal A}_{a, b, c}'=\sum_{i\in\z}\c v_i$ for some $a, b, c\in\c$, which is existed by the statements in \cite{LvZ} and Theorem 2.3. Then $V=U(H)V'$, where $U(H)=\c\{Y_{p_1}Y_{p_2}\cdots Y_{p_n}\mid \forall p_1, \cdots, p_n\in\z+s, \forall n\in\z_+\}$, as a vector space over $\c$.

\noindent{\bf Case 1. $c=0$}.  In this case $M_nV=0$. The irreducibility of $V$ as $\frak{sv}[s]$-modules is equivalent to that of $V$ as a $W(\frac12)[s]$-modules, where $W(\frac12)[s]=\c\{L_m, Y_{m+s}\mid m\in\z\}$ since $[Y_i, Y_j]v=(j-i)M_{i+j}v=0$. So by Theorem \ref{T41}, $V$ is the Harish-Chandra module of intermediate series.

\noindent {\bf Case 2. $c\ne 0$}.

Now set $U_i=U(H)v_i$ for any $i\in\z$.

\noindent{\bf Claim 1.} $U_i=U_j$ for all $i, j\in\z$. It means $V=U(H)v_i$ for any $i\in \z$.

In fact, for any $i\ne j\in\z$, $v_j=\frac1c M_{j-i}v_i=\frac1{c(j-i-2p)}[Y_p, Y_{j-i-p}]v_i\in U_i$ for some $p$ such that $j-i-2p\ne 0$, then $U_i\subset U_j$.  Similarly we have $U_j\subset U_i$. Then $U_i=U_j$.
So $V=U(H)v_i$ for any $i\in\z$. So we get the Claim 1. Moreover $M_j$ is nonsingular on $V$.

\noindent{\bf Claim 2.} All $Y_i$ are nonsingular on $V$.

In fact, if there exists $u_{i_0}\in V_i$ such that $Y_{p_0}u_{i_0}=0$ for some $p_0\in\z+s$. Then
we can easily to prove that $Y_{p_0}$ is nilpotent on $V_i$ and then it is nilpotent on $V$ since any $M_j$ is commutative with $Y_{p_0}$.
Suppose that $n$ is least positive integer such that $Y_{p_0}^nV=0$. By action of $Y_p\,(p\ne p_0)$ we have
$(p_0-p)Y_{p_0}^{n-1}M_{p+p_0}V=0$. Then $Y_{p_0}^{n-1}V=0$. It is contradict to the choice of $n$.

Now for any $h=Y_{p}^2M_{-2p}\in U(H)_0$ and $k\in(1-s)\z$, then $h$ is a linear transformation over the finite dimensional vector space $V_k$.
So there exists a minimal polynomial $m(x)$ of $h$ such that $m(h)V_k=0$ (if $s={1\over 2}$, then we can regard $h$ as a linear transformation over the finite dimensional vector space $V_k+V_{k+\frac12}$). By action some $M_j$ we have
\begin{equation}m(h)V=0.\label{nip2}\end{equation}

By action of $Y_{q}, (q\ne p)$, on (\ref{nip2}), we get $Y_{p}$ is not nonsingular on $V$. So we get a contradiction.\qed

\noindent{\bf Remark.}
Let $\mathcal S$ be the Lie subalgebra of $\frak{sv}[s]$ generated by $M_i, Y_p$ for all $i\in\z$ and $p\in\z+s$,
then $\mathcal S$ is just the infinite-dimensional Schr\"odinger algebra (see \cite{H}).
Clearly any irreducible uniformed bounded $\frak{sv}[s]$-module with all $Y_i$ nonsingular is just an irreducible uniform bounded $S$-module with all $Y_i$ nonsingular and $M_iu_j=cu_{i+j}$ for any $u_i\in V_i$. From the above proof in the Case 2, $\mathcal S$ has no such modules.

So from the above proposition, we have

\begin{theo}
Let $V$ be an irreducible weight ${\frak{sv}[s]}$-module with finite dimensional weight spaces. Then $V$ is a highest weight module, lowest weight module, or Harish-Chandra module of intermediate series.
\end{theo}

\subsection{The deformative twisted Schr\"{o}dinger-Virasoro Lie algebra}

For any $\varrho \in\mathbb Q$, \cite{RU} introduced a family of  infinite-dimensional Lie algebras called twisted deformative Schr\"{o}dinger-Virasoro Lie algebras ${\mathcal D}(\varrho)$,  admitting $\c$-basis
$\{L_n$,  $Y_n$,  $M_n\mid n\in\Z\}$ and the following Lie brackets
\begin{eqnarray*}
&&[L_m,L_n]=(n-m)L_{n+m},\\
&&[L_m,Y_n]=(n-\frac{\varrho+1}{2}m)Y_{n+m},\ \ \ [Y_n,Y_m]=(m-n)M_{n+m},\\
&&[L_m,M_n]=(n-\varrho m)M_{n+m},\ \ \ \ \ \ \, [Y_n,M_m]=[M_n,M_m]=0.
\end{eqnarray*}

Clearly, ${\mathcal D}(\varrho)$ contain $W(\varrho)$ as a subalgebra.
If $\varrho=0$, then ${\mathcal D}(\varrho)$ is just the above twisted Schr\"{o}dinger-Virasoro Lie algebra.
Now we suppose that $\varrho\ne0, -1, -3$ in this subsection.

\begin{prop}  Let $V$ be a uniformly bounded irreducible
${\mathcal D}(\varrho)$-module. Then $M_nV=Y_nV=0$. Then $V$ is the Harish-Chandra module of intermediate series.

\end{prop}

\noindent{\bf Proof.} Since $\varrho\ne 0$, the Lie subalgebra generated by $\{L_n, M_n\mid n\in\z\}$ is isomorphic to $W(\varrho)[0]$. Choose an irreducible $W(\varrho)[0]$-module $V'$ of $V$ and $M_nV'=0$ by Theorem 2.3.

Certainly as a ${\mathcal D}(\varrho)$-module $V=U(H)V'$, where $H=\c\{Y_n, n\in\z\}$. By $[M_n, Y_m]=0$ for all $m, n\in\z$, we have $M_nV=0$ and then all $[Y_m, Y_n]$ act as zero's since $[Y_n,Y_m]=(m-n)M_{n+m}$.
The irreducibility of $V$ as ${\mathcal D}(\varrho)$-modules is equivalent to that of $V$ as a $W(\frac{\varrho+1}2)[0]$-modules, where $W(\frac{\varrho+1}2)[0]=\c\{L_n, Y_n\mid n\in\z\}$ with $[Y_m, Y_n]=0$ for all $m, n\in\z$. So by Theorem \ref{T41}, $V$ is the Harish-Chandra module of intermediate series as an ${\mathcal D}(\varrho)$-module.\qed

Similarly we can also easily prove that an irreducible weight module $V$ is a uniform bounded module if it is not a highest weight module and a lowest weight module over ${\mathcal D}(\varrho)$. So we obtain the following result.

\begin{theo} Let $V$ be an irreducible weight ${\mathcal D}(\varrho)$-module with finite dimensional weight spaces. Then $V$ is a highest weight module, lowest weight module, or Harish-Chandra module of intermediate series.\qed
\end{theo}

\noindent{\bf Remarks.}
1. The Lie algebra ${\mathcal D}(1)$ is just the truncated Virasoro-loop algebra Vir$\ot$ $\c[t, t^{-1}]/(t^3)$. Its Harish-Chandra modules are classified in \cite{GLZ}.

2. Due to the central elements in the universal central extensions of the above Lie algebras are trivially act on the irreducible uniform bounded modules, so we have classified all Harish-Chandra modules over the universal central extensions of all the above Lie algebras.

3. With the methods and results in Section 3, we can classify Harish-Chandra modules over many Lie algebras, which consist of the Virasoro operators and some intertwining operators, such as Lie algebras listed in Section 4.

\vskip30pt \centerline{\bf ACKNOWLEDGMENTS}

\vskip15pt Project is supported by the NNSF (No. 11371134) and the
ZJNSF(LZ14A010001).

Part of the work was done by the author's visit in Department of Mathematics of Koeln University. He would like to express his special gratitude to
Prof. P. Littelmann for his kind hospitality and instructions, and to Dr. M. Fang for helpful discussions.

\def\refname{

\end{document}